\documentstyle[12pt]{article}
\author{Youssef ALAOUI}
\title{Holomorphic fiber bundle with Stein base and Stein fibers}
\date{}

\newcommand{\naturels}{\mbox{I}\!\!\!\mbox{N}}
\newcommand{\complexes}{\mbox{I}\!\!\!\mbox{C}}

\newtheorem{th}{theorem}
\newtheorem{lm}{lemma}
\newtheorem{crd}{Corollary}
\begin{document}
\maketitle
\setcounter{page}{1}
\noindent
{\large Summary}. In this article, we prove that
if $\Pi: X\rightarrow \Omega$ is a surjective holomorphic
map, with $\Omega$ a Stein space and $X$ a complex manifold
of dimension $n\geq 3,$ and if, for every $x\in \Omega$
there exists an open neighborhood $U$ such that $\Pi^{-1}(U)$
is Stein, then $X$ is Stein.\\ 
\\ 
{\large Key words}: Stein spaces; Stein mophism;
$q$-complete spaces; hyperconvex sets.\\ 
\\
$2000$ MS Classification numbers: 32E10, 32E40.\\
\\
\\
{\bf 1 Introduction}\\
\\
\hspace*{.1in}A surjective holomorphic map $\Pi: X\rightarrow Y$
between complex spaces is said to be a Stein morphism,
if for every point $y\in Y$, there exists an open neighborhood
$U$ of $y$ such that $\Pi^{-1}(U)$ is Stein.\\
\hspace*{.1in}In $1977$, Skoda $[12]$ showed that a locally
trivial analytic fiber bundle\\ $\Pi: X\rightarrow \Omega$ with
Stein base and Stein fibers is not necessarily a Stein manifold.
Thus giving a negative answer to a conjecture proposed in $1953$
by J-P. Serre $[11].$ 
In the counterexample provided by Skoda the base $\Omega$ is 
an open set in $\complexes$ and the holomorphic functions on 
$X$ are constant along the fibers 
$\Pi^{-1}(t)$, $t\in\Omega$. It follows (cf. J-P. Demailly $[4]$) 
that the cohomology group $H^{1}(X,{\mathcal{O}}_{X})$ is 
an infinitely dimensional complex vector space.\\
\hspace*{.1in}At the same time, Fornaess $[6]$ solved by means
of a $2$-dimensional counter-example the analogous problem for
the Stein morphism.\newpage
\noindent
\hspace*{.1in}There are still other counterexamples to the problem 
of Serre by J-P. Demailly $[5]$ in $1978$ and by Coeur\'e and Loeb $[3]$ in $1985$.\\
\hspace*{.1in}Note however that if $\Pi: X\rightarrow \Omega$
is a Stein morphism with a Stein base $\Omega$, it follows
from $[8]$ that for any
coherent analytic sheaf ${\mathcal{F}}$
on $X$ the cohomology groups $H^{p}(X,{\mathcal{F}})$ vanish 
for all $p\geq 2$. So it is therefore natural to raise the question
whether $X$ is $2$-complete.
The main new result of this article is to give
a positive answer to this problem when $X$ is a complex manifold.\\ 
\hspace*{.1in}We obtain this result as consequence of the following theorem
\begin{th}{-Let $X$ be a complex manifold of dimension $n\geq 3$, 
and $\Omega$ a complex space such that there exists a Stein
morphism $\Pi: X\rightarrow \Omega$. If $\Omega$ is Stein,
then $X$ itself is Stein.}
\end{th}
{\bf Proof of theorem $1$}\\
\\
\hspace*{.1in}We consider a covering ${\mathcal{V}}=(V_{i})_{i\in \naturels}$ of $\Omega$
by open sets $V_{i}\subset\Omega$ such that $\Pi^{-1}(V_{i})$
is Stein for all $i\in \naturels$. By the Stein covering lemma of
Sthel\'e $[13]$, there exits a locally finite covering 
${\mathcal{U}}=(U_{i})_{i\in \naturels}$ of $\Omega$ by Stein open subsets 
$U_{i}\subset\subset \Omega$ such that ${\mathcal{U}}$ is a refinement of ${\mathcal{V}}$ and
$\displaystyle\bigcup_{i\leq j}U_{i}$ is Stein for all $j$.
Moreover, there exists a continuous strictly psh function
$\phi_{j+1}$ in $\displaystyle\bigcup_{i\leq j+1}U_{i}$ such that
$$\displaystyle\bigcup_{i\leq j}U_{i}\cap U_{j+1}=\{x\in U_{j+1}: \phi_{j+1}(x)<0\}$$
Note also that $\Pi^{-1}(U_{i})$ is Stein for all $i\in \naturels$ and,
if $X_{j}=\Pi^{-1}(\displaystyle\bigcup_{i\leq j}U_{i})$ and
$X'_{j+1}=\Pi^{-1}(U_{j+1})$, then
$X_{j}\cap X'_{j+1}=\{x\in X'_{j+1}: \phi_{j+1}o\Pi(x)<0\}$ is Runge 
in $X'_{j+1}.$
\begin{lm}{-Under the above assumptions, the sets $X_{j}$, 
$j\in\naturels$, are Stein}
\end{lm} 
Proof. The proof is by induction on $j$.\\
\hspace*{.1in}For $j=0$, this is clear, since $\Pi^{-1}(U_{i})$
is Stein for all $i\in\naturels$. Now let $j\geq 1$ and
suppose that $X_{j}$ is Stein. Let 
$Y_{j}=\{x\in X_{j}: \phi_{j+1}o\Pi(x)>0\}$ and\\
$Y'_{j+1}=\{x\in X'_{j+1}: \phi_{j+1}o\Pi(x)>0\}.$ Then 
$H^{p}(Y_{j},{\mathcal{O}}_{X})\cong H^{p}(Y'_{j+1},{\mathcal{O}}_{X})=0$
for $1\leq p\leq n-2.$ In fact, 
let $\xi\in X'_{j+1}$ 
and $V\subset X'_{j+1}$ a hyperconvex open neighborhood of $\xi.$
Let $\psi: V\rightarrow ]-\infty,0[$ be a continuous strictly 
plurisubharmonic function. Then
$\psi_{k}=\frac{1}{k}\psi+\phi_{j+1}o \Pi,$ $k\geq 1,$ is an increasing
sequence of continuous strictly psh functions such that $\psi_{k}\rightarrow (\phi_{j+1}o \Pi)|_{V}.$
Let $V_{k}=\{x\in V: \psi_{k}(x)>0\},$ $k\geq 1.$
Then for every point $\xi'\in V$ such that $\psi_{k}(\xi')=0$
there exists, by [$2$, lemma $2$ ], a fundamental system
of connected Stein neighborhoods $U\subset\subset V$ of $\xi'$
such that
$H^{r}(U\cap V_{k}, {\mathcal{O}}_{X})=0$
for $1\leq r\leq n-2$ and the restriction map
$H^{o}(U, {\mathcal{O}}_{X})\rightarrow  H^{o}(U\cap V_{k}, {\mathcal{O}}_{X})$
is an isomorphism for every $k\in\naturels.$
Then a similar proof as in [$1$, lemma $2$] shows that
if $S'_{k}=\{x\in V: \psi_{k}(x)\leq 0\},$ then
$H^{p}_{S'_{k}}(V,{\mathcal{O}}_{X})=0$ for 
$p\leq n-1$, where $H^{p}_{S'_{k}}(V,{\mathcal{O}}_{X})$ is the 
p-th group of cohomology of $V$ with support in $S'_{k}$. Hence
the exact sequence of local cohomology
$$\cdots\rightarrow H^{p}_{S'_{k}}(V,{\mathcal{O}}_{X})
\rightarrow H^{p}(V,{\mathcal{O}}_{X})
\rightarrow H^{p}(V_{k},{\mathcal{O}}_{X})
\rightarrow H^{p+1}_{S'_{k}}(V,{\mathcal{O}}_{X})
\rightarrow\cdots$$
shows that $H^{p}(V_{k},{\mathcal{O}}_{X})=0$ for $1\leq p\leq n-2$
and\\ $H^{o}(V_{k},{\mathcal{O}}_{X})\cong H^{o}(V,{\mathcal{O}}_{X})$
for all $k\in \naturels.$
Since $V\cap Y'_{j+1}$ is an increasing union of $V_{k}$, $k\in \naturels$,
it follows from [$2$, lemma, p. $250$] that\\ $H^{p}(V\cap Y'_{j+1},{\mathcal{O}}_{X})=0$
for $1\leq p\leq n-2$ and 
$H^{o}(V, {\mathcal{O}}_{X})\rightarrow H^{o}(V\cap Y'_{j+1}, {\mathcal{O}}_{X})$
is an isomorphism. Since every point has a fundamental system
of hyperconvex neighborhoods, it follows from $[7]$ that
if\\ $S'_{j+1}=\{x\in X'_{j+1}: \phi_{j+1} o\Pi(x)\leq 0\},$
then $\underline{H^{p}_{S'_{j+1}}}({\mathcal{O}}_{X})=0$ for
$0\leq p\leq n-1,$ where  $\underline{H^{p}_{S'_{j+1}}}({\mathcal{O}}_{X})$
is the cohomology sheaf of $X'_{j+1}$ with coefficients in 
${\mathcal{O}}_{X}$ and support in $S'_{j+1}.$ Moreover,
there is a spectral sequence
$$H^{p}_{S'_{j+1}}(X'_{j+1}, {\mathcal{O}}_{X})\Longleftarrow
E_{2}^{p,q}=H^{p}(X'_{j+1}, \underline{H^{p}_{S'_{j+1}}}({\mathcal{O}}_{X}))$$
Since $\underline{H^{p}_{S'_{j+1}}}({\mathcal{O}}_{X})=0$ for
$p\leq n-1,$ then 
$H^{p}_{S'_{j+1}}(X'_{j+1}, {\mathcal{O}}_{X})=0$ for $p\leq n-1.$
Now by using the exact sequence of local cohomology
\begin{center}
$\cdots\rightarrow H^{p}_{S'_{j+1}}(X'_{j+1}, {\mathcal{O}}_{X}) 
\rightarrow H^{p}(X'_{j+1}, {\mathcal{O}}_{X}) 
\rightarrow H^{p}(Y'_{j+1}, {\mathcal{O}}_{X})
\rightarrow H^{p+1}_{S'_{j+1}}(X'_{j+1}, {\mathcal{O}}_{X})\rightarrow\cdots$ 
\end{center}
we see that 
$H^{p}(Y'_{j+1}, {\mathcal{O}}_{X})\cong H^{p}(X'_{j+1}, {\mathcal{O}}_{X})=0$
for $1\leq p\leq n-2.$\\
Similarly $H^{p}(Y_{j},{\mathcal{O}}_{X})=0$
for $1\leq p\leq n-2.$\\ 
We have 
$Y_{j+1}=\{x\in X_{j+1}: \phi_{j+1}o\Pi(x)>0\}=Y_{j}\cup Y'_{j+1}.$
Since $Y_{j}\cap Y'_{j+1}=\emptyset$, then 
$H^{p}(Y_{j+1},{\mathcal{O}}_{X})=0$ for $1\leq p\leq n-2.$\\
On the other hand, if $S_{j+1}=\{x\in X_{j+1}: \phi_{j+1} o\Pi(x)\leq 0\},$
then we can prove exactly as for $H^{p}_{S'_{j+1}}(X'_{j+1}, {\mathcal{O}}_{X})$ 
that $H^{p}_{S_{j+1}}(X_{j+1}, {\mathcal{O}}_{X})=0$ for $p\leq n-1.$ 
Now the exact sequence of local cohomology
$$\cdots\rightarrow H^{1}_{S_{j+1}}(X_{j+1},{\mathcal{O}}_{X})    
\rightarrow H^{1}(X_{j+1},{\mathcal{O}}_{X})\rightarrow H^{1}(Y_{j+1},{\mathcal{O}}_{X})=0$$
implies that $H^{1}(X_{j+1},{\mathcal{O}}_{X})=0.$ 
This proves that $X_{j+1}$ is Stein. (See [9]).\\
\\
{\bf End of the proof of theorem $1$}\\
\\
\hspace*{.1in} We consider for every $j\geq 1$ the long exact sequence 
of cohomology associated to the Mayer-Vietoris sequence
\begin{center}
$\cdots\rightarrow H^{o}(X_{j+1},{\mathcal{O}}_{X})\rightarrow
H^{o}(X_{j},{\mathcal{O}}_{X})\oplus H^{o}(X'_{j+1},{\mathcal{O}}_{X})
\rightarrow H^{o}(X_{j}\cap X'_{j+1},{\mathcal{O}}_{X})
\rightarrow H^{1}(X_{j+1},{\mathcal{O}}_{X})\rightarrow\cdots$
\end{center}

Since $H^{1}(X_{j+1},{\mathcal{O}}_{X})=0$ and 
the restriction map
$$H^{o}(X'_{j+1},{\mathcal{O}}_{X})\rightarrow H^{o}(X_{j}\cap X'_{j+1},{\mathcal{O}}_{X})$$
has dense image, then, by ($[2]$, proof of Proposition $19$),
one deduces that $X_{j}$ is Runge in $X_{j+1},$ which implies (see $[14]$) 
that $X=\displaystyle\bigcup_{j\geq 1}X_{j}$ is Stein.\\
\\ 
\begin{crd}{-Let $\Pi: X\rightarrow \Omega$ be a Stein morphism
with a Stein base $\Omega.$
Assume that\\
\hspace*{.1in}a) $X$ is a complex manifold of dimension $n\geq 3$\\
or\\
\hspace*{.1in}b) $X$ a complex space of dimension $2.$\\
Then $X$ is $2$-complete.}
\end{crd}
{\bf Proof}\\
\\ 
\hspace*{.1in}If $dim(X)=n\geq 3,$ then, by theorem $1,$ $X$ is Stein.\\
\hspace*{.1in}Suppose that $dim(X)=2$ and let $A$ be a compact connected analytic
subset of $X.$  Since $\Omega$ is Stein, then $\Pi(A)$ is reduced
to one point. This implies that $A$ is contained in one fiber.
Since the fibers of $\Pi$ are Stein, then $A$ is $0$-dimensional,
which proves by a theorem of Oshawa $[10]$ that $X$ is $2$-complete.


\begin{thebibliography}{10}

\bibitem{bib1} Y. Alaoui, Cohomology of locally $q$-complete sets in
Stein manifolds. Complex Variables and Elliptic Equations.
Vol. $51,$ No. $2,$ February $2006,$ $137-141$\\
\\
\bibitem{bib2} A. Andreotti and H. Grauert, 
Th\'eor\`emes de finitude de la cohomologie des espaces
complexes. Bull. Soc. Math. France $90$ ($1962,$) $193-259$\\ 
\\
\bibitem{bib3}G. Coeur\'e and J. J. Loeb, A counterexample
to the Serre problem with a bounded domain of $\complexes^{2}$
as fiber. Annals of Mathematics, $122$ ($1985,$) $329-334$\\
\\
\bibitem{bib4} J-P. Demailly, Diff\'erents examples de fibr\'es 
holomorphes non de Stein. S\'eminaire de Lelong-Skoda,
Lecture Notes in Math num\'ero $694,$ $1976-77,$ p. $15-41$\\
\\
\bibitem{bib5} J-P. Demailly, Un exemple de fibr\'e holomorphe
non de Stein \`a fibre $\complexes^{2}$ ayant pour base le disque
ou le plan. Inventiones mathematicae $48,$ $293-302$ ($1978$)\\ 
\\
\bibitem{bib6}J-E. Fornaess, $2$-Dimensional Counterxamples to
Generalizations of the Levi Problem, 
Mathematische Annalen $230,$ $169-174$ ($1977$)\\
\\
\bibitem{bib7}A. Grothendieck, $1957,$ Sur quelques points
d'Alg\`ebre homologique. Tohuku Mathematical Journal, IX,
$119-221.$\\
\\
\bibitem{bib8}B. Jennane, Groupes de cohomologie d'un fibr\'e
holomorphe \`a base et \`a fibre de Stein,
Inventiones Mathematicae $54,$ $75-79$ $(1979)$\\
\\
\bibitem{bib9}B. Jennane, Probl\`eme de Levi et morphisme
localement de Stein,
Math. Ann. $256,$ $37-42$ $(1981)$\\
\\
\bibitem{bib10}T. Ohsawa, Completeness of non-compact analytic 
spaces, Publ. R. I. M. S. Kyoto Univ. $20$ ($1984$), $683-692.$\\
\\
\bibitem{bib11}J-P. Serre, Quelques probl\`emes  globaux relatifs
aux vari\'et\'es de Stein. Colloque sur les fonctions de
plusieurs variables, Bruxelles, $1953.$\\
\\
\bibitem{bib12}H. Skoda, Fibr\'es holomorphes \`a base et \`a fibre
de Stein. Inventiones Mathematicae $43,$ $97-107$ $(1977)$\\
\\
\bibitem{bib13}J-L. Stehl\'e, Fonctions plurisousharmoniques
et conv\'exit\'e holomorphe de certains fibr\'es analytiques.
S\'eminaire Lelong, p. $155-179,$ Lecture Notes num\'ero
$474,$ $1973,74$\\
\\
\bibitem{bib14}K, Stein, \"Uberlagerungen holomorph-vollst\"andiger komplexer
R\"aume. Arch. Math. $7,$ $356-361$ ($1956$)
\end{thebibliography}
\end{document}